\documentclass[11pt]{article}
\usepackage{newlfont,amsfonts,amssymb,oldgerm,amsmath,amsthm,amsgen,amscd}
\begin{document}

\newcommand{\Q}{\ensuremath{\mathbb{H}}}
\newcommand{\N}{\ensuremath{\mathbb{N}}}
\newcommand{\Z}{\ensuremath{\mathbb{Z}}}
\newcommand{\C}{\ensuremath{\mathbb{C}}}
\newcommand{\K}{\ensuremath{\mathbb{K}}}
\renewcommand{\O}{\ensuremath{\mathcal{O}}}
\newcommand{\R}{\ensuremath{\mathbb{R}}}

\newcommand{\bcase}{\begin{case}}
\newcommand{\ecase}{\end{case}}
\newcommand{\setcase}{\setcounter{case}{0}}
\newcommand{\bclaim}{\begin{claim}}
\newcommand{\eclaim}{\end{claim}}
\newcommand{\setclaim}{\setcounter{claim}{0}}
\newcommand{\bstep}{\begin{step}}
\newcommand{\estep}{\end{step}}
\newcommand{\setstep}{\setcounter{step}{0}}
\newcommand{\bhlem}{\begin{hlem}}
\newcommand{\ehlem}{\end{hlem}}
\newcommand{\sethlem}{\setcounter{hlem}{0}}

\newcommand{\bleer}{\begin{leer}}
\newcommand{\eleer}{\end{leer}}
\newcommand{\bde}{\begin{de}}
\newcommand{\ede}{\end{de}}
\newcommand{\ul}{\underline}
\newcommand{\ol}{\overline}
\newcommand{\tbf}{\textbf}
\newcommand{\mc}{\mathcal}
\newcommand{\mb}{\mathbb}
\newcommand{\mf}{\mathfrak}
\newcommand{\bs}{\begin{satz}}
\newcommand{\es}{\end{satz}}
\newcommand{\btheo}{\begin{theo}}
\newcommand{\etheo}{\end{theo}}
\newcommand{\bfolg}{\begin{folg}}
\newcommand{\efolg}{\end{folg}}
\newcommand{\blem}{\begin{lem}}
\newcommand{\elem}{\end{lem}}
\newcommand{\bnote}{\begin{note}}
\newcommand{\enote}{\end{note}}
\newcommand{\bprf}{\begin{proof}}
\newcommand{\eprf}{\end{proof}}
\newcommand{\bd}{\begin{displaymath}}
\newcommand{\ed}{\end{displaymath}}
\newcommand{\be}{\begin{eqnarray*}}
\newcommand{\ee}{\end{eqnarray*}}
\newcommand{\eeqa}{\end{eqnarray}}
\newcommand{\beqa}{\begin{eqnarray}}
\newcommand{\bi}{\begin{itemize}}
\newcommand{\ei}{\end{itemize}}
\newcommand{\bnum}{\begin{enumerate}}
\newcommand{\enum}{\end{enumerate}}
\newcommand{\la}{\langle}
\newcommand{\ra}{\rangle}
\newcommand{\eps}{\epsilon}
\newcommand{\ve}{\varepsilon}
\newcommand{\vp}{\varphi}
\newcommand{\lra}{\longrightarrow}
\newcommand{\Lra}{\Leftrightarrow}
\newcommand{\Ra}{\Rightarrow}
\newcommand{\sub}{\subset}
\newcommand{\ems}{\emptyset}
\newcommand{\sms}{\setminus}
\newcommand{\ints}{\int\limits}
\newcommand{\sums}{\sum\limits}
\newcommand{\lims}{\lim\limits}
\newcommand{\bcup}{\bigcup\limits}
\newcommand{\bcap}{\bigcap\limits}
\newcommand{\beq}{\begin{equation}}
\newcommand{\eeq}{\end{equation}}
\newcommand{\einhalb}{\frac{1}{2}}
\newcommand{\rr}{\mathbb{R}}
\newcommand{\rn}{\mathbb{R}^n}
\newcommand{\ccc}{\mathbb{C}}
\newcommand{\cn}{\mathbb{C}^n}
\newcommand{\M}{{\cal M}}
\newcommand{\drehgleich}{\mbox{\begin{rotate}{90}$=$  \end{rotate}}}
\newcommand{\turngleich}{\mbox{\begin{turn}{90}$=$  \end{turn}}}
\newcommand{\turnsimeq}{\mbox{\begin{turn}{270}$\simeq$  \end{turn}}}
\newcommand{\vf}{\varphi}
\newcommand{\earr}{\end{array}\]}
\newcommand{\barr}{\[\begin{array}}
\newcommand{\bvec}{\left(\begin{array}{c}}
\newcommand{\evec}{\end{array}\right)}
\newcommand{\sumk}{\sum_{k=1}^n}
\newcommand{\sumi}{\sum_{i=1}^n}
\newcommand{\suml}{\sum_{l=1}^n}
\newcommand{\sumj}{\sum_{j=1}^n}
\newcommand{\sumij}{\sum_{i,j=1}^n}
\newcommand{\suminf}{\sum_{k=0}^\infty}
\newcommand{\inv}{\frac{1}}
\newcommand{\wzbw}{\hfill $\Box$\\[0.2cm]}
\newcommand{\lag}{\mathfrak{g}}
\newcommand{\lan}{\mathfrak{n}}
\newcommand{\lah}{\mathfrak{h}}
\newcommand{\laz}{\mathfrak{z}}
\newcommand{\+}{\oplus}
\newcommand{\x}{\times}
\newcommand{\lx}{\ltimes}
\newcommand{\rrn}{\mathbb{R}^n}
\newcommand{\laso}{\mathfrak{so}}
\newcommand{\lason}{\mathfrak{so}(n)}
\newcommand{\lagl}{\mathfrak{gl}}
\newcommand{\lasl}{\mathfrak{sl}}
\newcommand{\lasp}{\mathfrak{sp}}
\newcommand{\lasu}{\mathfrak{su}}
\newcommand{\w}{\omega}
\newcommand{\pmh}{{\cal P}(M,h)}
\newcommand{\s}{\sigma}
\newcommand{\deri}{\frac{\partial}}
\newcommand{\ddx}{\frac{\partial}{\partial x}}
\newcommand{\ddz}{\frac{\partial}{\partial z}}
\newcommand{\ddi}{\frac{\partial}{\partial y_i}}
\newcommand{\ddj}{\frac{\partial}{\partial y_j}}
\newcommand{\ddk}{\frac{\partial}{\partial y_k}}
\newcommand{\ddp}{\frac{\partial}{\partial p_i}}
\newcommand{\ddq}{\frac{\partial}{\partial q_i}}
\newcommand{\xz}{^{(x,z)}}
\newcommand{\mh}{(M,h)}
\newcommand{\wxz}{W_{(x,z)}}
\newcommand{\qmh}{{\cal Q}(M,h)}
\newcommand{\bbem}{\begin{bem}}
\newcommand{\ebem}{\end{bem}}
\newcommand{\bbez}{\begin{bez}}
\newcommand{\ebez}{\end{bez}}
\newcommand{\bbsp}{\begin{bsp}}
\newcommand{\ebsp}{\end{bsp}}
\newcommand{\pr}{pr_{\lason}}
\newcommand{\huts}{\hat{\s}}
\newcommand{\whut}{\w^{\huts}}
\newcommand{\bhg}{{\cal B}_H(\lag)}
\newcommand{\aaa}{\alpha}
\newcommand{\bb}{\beta}
\newcommand{\laa}{\mf{a}}
\newcommand{\lam}{\lambda}
\newcommand{\LL}{\Lambda}
\newcommand{\D}{\Delta}
\newcommand{\ß}{\beta}
\newcommand{\ä}{\alpha}
\newcommand{\W}{\Omega}
\newcommand{\esel}{\ensuremath{\mathfrak{sl}(2,\ccc)}}
\newcommand{\kg}{{\cal K}(\lag)}
\newcommand{\bg}{{\cal B}_h(\lag)}
\newcommand{\kk}{  \mathbb{K}}
\newcommand{\xy}{[x,y]}
\newcommand{\perdef}{$\stackrel{\text{\tiny def}}{\iff}$}
\newcommand{\eqdef}{\stackrel{\text{\tiny def}}{=}}
\newcommand{\lai}{\mf{i}}
\newcommand{\lar}{\mf{r}}
\newcommand{\Dim}{\mathsf{dim\ }}
\newcommand{\im}{\mathsf{im\ }}
\newcommand{\Ker}{\mathsf{ker\ }}
\newcommand{\trace}{\mathsf{trace\ }}
\newcommand{\grad}{\mathsf{grad}}
\newcommand{\wt}{\widetilde}
\newcommand{\tnab}{\widetilde{\nabla}}
\newcommand{\tem}{\widetilde{M}}
\newcommand{\nabt}{\nabla^{\cal T}}
\newcommand{\ro}{\mathsf{P}}
\newcommand{\lecturecount}{\begin{center}{\sf  (Lecture \Roman{lecturenr})}\end{center}\addtocontents{toc}{{\sf  (Lecture \Roman{lecturenr})}} \refstepcounter{lecturenr}
}
%% CHOL
\newcommand{\T}{{\cal T}}
\newcommand{\cur}{{\cal R}} 
\newcommand{\pd}{{\cal P}} 

\newcommand{\inter}{\makebox[11pt]{\rule{6pt}{.3pt}\rule{.3pt}{5pt}}}

\newcommand{\onabla}{\overrightarrow{\nabla}}
\newcommand{\oR}{\overrightarrow{R}}
\newcommand{\Ric}{\mathsf{Ric}}

%\newcommand{\insrt}{\mbox{\begin{turn}{180} $\invneg$  \end{turn}}}
%\newcommand{\turngleich}{\mbox{\begin{turn}{90}$=$  \end{turn}}}
%\newcommand{\turnsimeq}{\mbox{\begin{turn}{270}$\simeq$  \end{turn}}}

%\swapnumbers
\theoremstyle{definition}
\newtheorem{de}{Definition}[section]
\newtheorem{bem}[de]{Remark}
\newtheorem{bez}[de]{Notation}
\newtheorem{bsp}[de]{Example}
\theoremstyle{plain}
\newtheorem{lem}[de]{Lemma}
\newtheorem{satz}[de]{Proposition}
\newtheorem{folg}[de]{Corollary}
\newtheorem{theo}[de]{Theorem}

\bibliographystyle{alpha}

%\markright{\centerline{{\sf --- Preliminary version ---}}}

\title{Ambient connections realising conformal Tractor holonomy}

\author{Stuart Armstrong\thanks{stuart.armstrong@st-cross.oxford.ac.uk } \ and Thomas Leistner\thanks{
%{\sc School of Mathematical Sciences, The University of Adelaide, SA 5006, Australia.} {\em email:} 
tleistne@maths.adelaide.edu.au
\newline {\em Date:} June 15, 2006 
}}
%\\[.3cm]
%{\small {\em Pure Mathematics - School of Mathematical Sciences}}\\
%{\small {\em The University of Adelaide, SA 5005, Australia}}\\[.2cm]
%{\small {\em phone: +61 (0)8 83033712,
%fax: +61 (0)8 83033696}}\\
%{\small {\em email:} {\tt tleistne@maths.adelaide.edu.au}}
%}
\date{ }
\maketitle
 \begin{abstract}
For a conformal manifold we introduce the notion of an ambient connection, an affine connection on an ambient manifold of the conformal manifold, possibly with torsion, and with conditions relating it to the conformal structure. The purpose of this construction is to realise the normal conformal tractor holonomy as affine holonomy of such a connection. We give an example of an ambient connection for which this is the case, and which is torsion free if we start the construction with a C-space, and in addition Ricci-flat if we start with an Einstein manifold. Thus for a $C$-space this example leads to an ambient metric in the weaker sense of \v{C}ap and Gover, and for an Einstein space to a Ricci-flat ambient metric in the  sense of Fefferman and Graham.
\end{abstract}

%\tableofcontents

\section{Introduction}

Conformal geometry has been studied for a long time. Cartan's \cite{ECC} techniques of `moving frames' developed into principal bundles and the Cartan connection.
The equivalent `Tractor' connection has been developed by T. Thomas (\cite{thomas26}   and \cite{thomas32}), further developed by T.N. Bailey, M.G. Eastwood and A.R. Gover (\cite{bailey-eastwood-gover94} and \cite{eastwood95}), and extensively treated in papers of A.R Gover and A. \v{C}ap (e.g. \cite{cap/gover02} and \cite{cap-gover03}), which being a vector bundle construction, allows for more explicit calculations. One interesting invariant of the Tractor connection is its holonomy group, and the geometric consequences of this group. Classification of these groups is essentially finished for 
Riemannian  signature in \cite{armstrong05}, see also \cite{leitner04killing}.

A major technique in this classification was the use of cone-constructions for conformally Einstein manifolds with non-zero Einstein constant, and a degenerate cone construction in the Ricci-flat case \cite{leistner05}, which are both instances of  the original ambient construction of C. Fefferman and C. R. Graham (\cite{fefferman/graham85}, 
 also \cite{fefferman/graham02} and \cite{fefferman/hirachi03}). A qualitative description generalising the original ambient construction was given in \cite{cap-gover03}, where also the relation to the Tractor bundle is described.
This paper seeks to generalise both cone-constructions and the ambient metric construction of \cite{cap-gover03} to generate an affine connection on an ambient manifold which is metric but \emph{possibly with torsion} and  with same holonomy group as the Tractor connection -- indeed the tangent bundle of this ambient manifold becomes identified with the Tractor bundle.

In this paper, we will start by introducing the basic concepts of 
tractor bundles in Section 2. We shall also define associated structures, and 
set up the notation we shall be using throughout.
Section 3 defines recalls the notion of an ambient manifold with ambient metric following \cite{cap-gover03}, introduces the notion of an {\em ambient connection}, and demonstrates many of the properties it must have. Section 4 then looks more particularly at the intermesh between an ambient manifold and the Tractor bundle and how the 
properties of one transfer across to the other.

What is presented in Section 5 was the inspiration for this paper: an example of an ambient connection is given the holonomy of which equals to the normal conformal Tractor holonomy, and that, we hope, could lead to useful results for calculating Tractor holonomy groups. 
For now it serves two purposes: on the one hand
it demonstrates geometrically the algebraic fact proved in \cite{leistner05} that the Tractor holonomy algebra of a conformal structure which contains a metric with vanishing Cotton--York tensor is a Berger algebra. On the other hand our construction contains the Einstein cone and Ricci-flat construction as special cases.
Finally, in Section 6 we give some useful properties of the holonomy of ambient connections, for future constructions.

\paragraph{Acknowledgements.} The authors acknowledge the hospitality at the Erwin Sch\"{o}dinger International Institute for Mathematical Physics (ESI) in Vienna, and would like to thank the  organisers of the programme {\em Geometry of Pseudo-Riemannian Manifolds with Applications in Physics} at the ESI, during which the  discussion on the topic started.

 \section{Preliminaries}

Let $M$ be a smooth manifold of dimension $n$.

\bde[Weight Bundles] Weight bundles are defined as line bundles that are powers of the top wedge product of the cotangent bundle. In details,
\be
\cal E  [-n] = \wedge^n T^*M
\ee
and
\be
\cal E  [a] = \big( \cal E \cal [-n] \big)^{-\frac{a}{n}}.
\ee
\ede

If $\cal B $ is a vector bundle associated to the tangent bundle, we designate the tensor product $\cal B \otimes \cal E [a]$ as $\cal B \cal [a]$. Weight bundles are a way of formalising the notion of conformal weights in the following sense.
If we are given a conformal structure $[g]$ as a equivalence class of metrics with $\rr^+$ action
\be
(t, g) \to t^2 g,
\ee
then the action on the determinant becomes
\be
(t, \mathsf{det}( g)) \to t^n \mathsf{det} (g)
\ee
and since $\mathsf{det}( g) \in \Gamma(\cal E [-n])$, it can be seen that the conformal structure $[g]$ is isomorphic to a non-degenerate section $\hat{g}$ of $(\odot^2 T^*M)[2]$.

%For a manifold with conformal structure the frame bundle reduces to the bundle of frames orthonormal with respect to a metric in the conformal class, denoted by $\cal G_0$. The conformal group  $G_0:=CO(p,q)$ is the structure group of this bundle. A linear connection of this bundle is called {\em preferred} if has no torsion. 
There are many equivalent ways of describing the conformal structure, quite
apart from  equivalence classes of proportionate metrics, such as classes of
conformal connections, a principal $CO(p,q)$-bundle for the tangent bundle, etc.
We will use the description in terms of the  {\em standard Tractor bundle} (see \cite{thomas26}, \cite{bailey-eastwood-gover94}, and \cite{cap/gover02}) equipped with the standard tractor connection.
This construction defines the conformal structure entirely; however, for a given
conformal structure, the Tractor connection is not unique. Similarly to the
Levi-Civita connection for metric structures, there is a unique
{\em normal}  Tractor
connection for each conformal structure. For our purposes a conformal standard Tractor bundle with Tractor connection and the normal connection amongst them are best described by the following properties, from which the conformal structure can be recovered.

\begin{theo}[\cite{cap-gover03},\cite{cap/gover02}]
\label{capgoverthm}
Let $M$ be a smooth  manifold of dimension $n=p+q\ge 3$  and $\cal T$ be a vector bundle over $M$ of rank $n+2$.
\bnum
\item
 Suppose $\cal T$ admits a metric $h$ of signature $(p+1,q+1)$, an injective vector bundle homomorphism $\cal E[-1]\hookrightarrow\cal T$ with image $\cal T^1$, and a covariant derivative $\onabla$ such that
\bnum
\item $\cal T^1$ is light-like,
\item $\onabla h=0$, and
\item the following non-degeneracy condition is satisfied: for any $x\in M$ and any section $\s \in \Gamma(\cal T^1)$, non-zero at $x$, the map
$
\onabla \s : T_x M\to \cal T_x
$
is non-degenerate.
\enum
Then $(\cal T^1)^\bot / \cal T^1\simeq TM[-1]$ and $(\cal T, h, \onabla)$ is a standard tractor bundle for the conformal structure defined by the restriction of $h$ onto $\left(\cal T^1\right)^\bot / \cal T^1$.

\item A standard conformal Tractor connection $\onabla$ on $\cal T$ is {\em normal}  if and only if its curvature $\cal R$ maps $\cal T^1$ onto $ \cal T^1$ and the Ricci contraction of $W\in \Lambda^2 T^*M \otimes End( TM[-1])$, which is defined by $W(X,Y)[Z]=[\cal R(X,Y)Z]$ using the identification $TM[-1]\simeq \left(\cal T^1\right)^\bot / \cal T^1$, vanishes. With these normality conditions $(\cal T,\cal T^1, h, \onabla)$ is uniquely determined by the conformal structure up to isomorphism.
\enum
\etheo
To see the first point, notice that the conditions imply that $\onabla \s \in \Gamma(TM\otimes (\cal T^1)^\bot)$ and that for any choice of $\s$, one has a splitting of
\beqa \label{T:splitting}
\cal T = \cal E [1] \oplus TM[-1] \oplus \cal E[-1],
\eeqa
where $TM[-1]$ is the image of $TM$ under $\onabla \s$, and $\cal E[1]$ is the span of a section $\s'$, with $\s' \bot TM[-1]$ and $h(\s', \s) = -1$. This allows us to project $\onabla$ to get a connection $\nabla$ on $TM[-1]$, which conserves $\hat{g} = h|_{TM[-1]}$ --- in other words, conserves a conformal structure. Since $TM[-1] \oplus E[-1] \cong (\cal T^1)^{\bot}$, different non-zero choices of $\s$ give the same $\hat{g}$ --- though different $\nabla$'s.

Using $\s$ to identify $TM[-1]$ with $TM$, $\nabla$ then descends to a tangent bundle connection, preserving the tangent metric $\s^{-2} \hat{g}$. We call these $\nabla$ preferred connections of the Tractor connection $\onabla$.
The invariance of $\cal T^1$ under the curvature of the Tractor connection then
is equivalent with the preferred connections $\nabla$ being
torsion-free. This is because, in the previous notation, the  $TM[-1] \cong (\cal T^1)^\bot/\cal T^1$ component of $R(X,Y) \s$ is just $\nabla_X Y - \nabla_Y X - [X,Y]$.
%
%(\sf{I think we don't need the Ricci curvature condition}).
%
%and the Ricci contraction of $W\in \Lambda^2 T^*M \otimes End( TM[-1])$, which is defined by $W(X,Y)[Z]=[\cal R(X,Y)Z]$ using the identification $TM[-1]\simeq \left(\cal T^1\right)^\bot / \cal T^1$ vanishes. 
%\ede

So now pick a preferred $\nabla$, which preserves a metric $g$, and (since it is torsion-free) is the Levi-Civita connection of $g$. Given the associated splitting
\be
\cal T = \cal E [1] \oplus TM[-1] \oplus \cal E[-1],
\ee
the normal Tractor connection $\onabla$ is given as
\beqa \label{Tractor:formula}
\onabla_X = \nabla_X + X + \mathsf{P}(X).
\eeqa
Here $\mathsf{P}$ is the \emph{Schouten} tensor of $\nabla$ defined by
\begin{eqnarray*}
\mathsf{P}_{ij} = - \frac{1}{n-2} \ \big( \Ric_{ij} \ - \ \frac{1}{2n-2} R \hat{g}_{ij} \big),
\end{eqnarray*}
$\Ric$ being the Ricci tensor and $R$ the scalar curvature of $\nabla$. To understand the action of $X$ and $\mathsf{P}(X)$ in (\ref{Tractor:formula}), let $\nu, A$ and $\s$ be sections of $\cal E[1]$, $TM[-1]$ and $\cal E [-1]$, respectively. Then the action is
\be
\begin{array}{lcllclcll}
X \cdot \nu & = & 0 && \ \ \ & \mathsf{P}(X) \cdot \nu &=& \nu \hat{g}^{-1}(\mathsf{P}(X)) &\in \Gamma(TM[-1]) \\
X \cdot A & = & \hat{g}(X,A) &\in \Gamma(\cal E[1]) & \ \ \ & \mathsf{P}(X) \cdot A &=& \mathsf{P}(X) \inter A &\in \Gamma(\cal E[-1]) \\
X \cdot \s & = & \s X &\in \Gamma(TM[-1]) & \ \ \ & \mathsf{P}(X) \cdot \s &=& 0.&
\end{array}
\ee

\section{Ambient manifolds, ambient metrics, and ambient connections}
\label{section-ambient}
Let $ (M,c)$ be a  a smooth $n$-dimensional manifold $M$ with conformal structure $c$ of signature $(p,q)$. This is equivalent to a reduction of the linear frame bundle $Gl(M)$ to the conformal frame bundle $CO(M)$. It can also be characterised by a principle $\rr^+$-fibre bundle $\pi: \cal Q\rightarrow M$ defined as the ray sub-bundle in the bundle of metrics of signature $(p,q)$ given by metrics in the conformal class $c$. In fact, via the $\mathsf{det}^{w/n}$--reduction of $Gl(M)$ or $CO(M)$ to $\cal Q$, it can be viewed as the principal bundle to the weight bundles defined above. The action of $\rr^+$ on $\cal Q$ shall be denoted by $\vf$:
\be
\vf(t,g_x)&=& t^2 g_x.\ee 
From \cite{cap-gover03} we recall the following definitions.
\bde\label{ambientmet}
Let $\pi:\cal Q\rightarrow M$ be a conformal structure of signature $(p,q)$ over an $n$-dimensional manifold $M$. 
\bnum
\item An $(n+2)$-dimensional manifold $\widetilde{M}$ is called {\em ambient manifold} if
\bnum
\item there is a free $\rr^+$-action $\widetilde{\vf}$ on $M$, and
\item an embedding $\iota:\cal Q\rightarrow \wt{M}$ which is $\rr^+$-equivariant, i.e. the diagram
\[ \begin{CD}
\rr^+ \times \cal Q @>{\vf}>> \cal Q \\
@V{id\times\iota}VV  @VV{\iota}V\\
\rr^+\times \wt{M} @>{\wt{\vf}}>> \wt{M}
\end{CD}\]
commutes.
\enum
\item A metric $h$ on an ambient manifold $\wt{M}$ with $\rr^+$-action $\wt{\vf}$ is called {\em ambient metric}, if the following is satisfied:
\bnum
\item If $F$ is the fundamental vector field of $\wt{\vf}$, and $\cal L$ denotes the Lie derivative, then
\beqa\label{methom}
\cal L_Fh&=&2 h,\eeqa
i.e. the metric is homogeneous of degree $2$ w.r.t. the $\rr^+$-action.
\item For any $g_x\in \cal Q$ the following equality in $\odot^2 T^*_{g_x}\cal Q$ holds:
\beqa\label{metext}
(\iota^*h)_{g_x}& =& g_x \left(d\pi(.),d\pi(.)\right).\eeqa
\enum
\enum
\ede
We should point out that in the original definition of an ambient metric by C. Fefferman and C.R. Graham in \cite{fefferman/graham85}, the additional condition of Ricci-flatness is imposed.

Recall that for $x \in \wt{M}$, the fundamental vector field $F(x)$ of the $\rr^+$-action is defined as the tangent vector of the curve $ \wt{\vf}(\text{e}^{t},x)$ at $t=0$. Its flow through $x\in \tem$ is given by $\vf_{e^t}(x)$.
Let us define the $1$-form $\phi$ as the dual to the fundamental vector field $F$, i.e. $\phi:=h(F,.)$. Since $F$ is the image under $\iota^*$ of the fundamental vector field of the action $\vf$ on $\cal Q$, (\ref{metext}) gives that 
\beqa \iota^*\phi &=&0.
\eeqa
This implies that $F$ is of length zero along $\iota(\cal Q)$ and that $h$ is of signature $(p+1,q+1)$. Furthermore it implies that
\beqa d\iota^*\phi\ =\ \iota^* d\phi&=&0.\eeqa
Next, we want to consider sections which are homogeneous w.r.t.~the $\rr^+$ action on $\tem$. We introduce the following notation: If $\cal E$ is a tensor bundle over the ambient manifold $\tem$ or over $\iota(\cal Q)$, then we set for an integer $w$
\[ \Gamma^{w}(\cal E)\ := \ \left\{ s\in \Gamma(\cal E) \mid \cal L_F s= ws\right\}.\]
We write $\Gamma^w(\rr)$ for functions which are homogeneous of degree $w$, i.e.
$f\in \Gamma^w(\rr)$ if $\cal L_F f=F(f)= w\cdot f $. 
Direct calculation using the definition of the Lie derivative and the Jacobi identity for the commutator gives the following properties.

\blem \label{tensordegrees}
Let $\tem$ be a manifold with $\rr^+$ action and corresponding fundamental vector field $F$. Then it holds:
\bnum
\item
If $\psi\in \Gamma^w(\otimes^kT^*\tem)$, and $X_i\in \Gamma^{w_i}(T\tem)$ for $i=1, \ldots , k$  the function $\psi(X_1, \ldots  , X_k)$ is homogeneous of degree $w+w_1+\ldots + w_k$.
\item If $X\in \Gamma^v(T\tem)$ and $Y\in \Gamma^w(T\tem)$, then $\left[X,Y\right]\in \Gamma^{v+w}(\T\tem)$.
\item If $f\in \Gamma^v(\rr)$ and $X\in \Gamma^w(T\tem)$, then $X(f)\in \Gamma^{v+w}(\rr)$.
\enum
\elem
\bigskip

Of course, the ambient metric defines uniquely the Levi-Civita connection, but for our purposes it is more useful to work with a more general connection which still parallelises the ambient metric.
\bde
Let $(\wt{M}, h)$ be an ambient manifold with ambient metric $h$ of a conformal structure $\cal Q\rightarrow M$ in the sense of Definition \ref{ambientmet}. A linear connection  $\wt{\nabla}$ is called {\em ambient connection} if
$\tnab h\ =\ 0$
and the $(3,0)$-torsion $T^*$ of $\wt{\nabla}$ satisfies
\beqa\cal L_FT^*\ =\ 2 T^*,\label{conhom}\eeqa
i.e. is homogeneous of degree $2$ w.r.t. the $\rr^+$-action on $\wt{M}$.
\ede
Recall that the $(2,1)$-torsion of $\wt{\nabla}$ is defined as 
$T_{X,Y}:=\tnab_XY-\tnab_YX- [X,Y]$
and the $(3,0)$-torsion as $T^*(X,Y,Z):=h(T_{X,Y},Z)$. Naturally,  (\ref{conhom}) is equivalent to $\cal L_F T=0$.

If we define the $1$-form $\phi$ as the dual to the fundamental vector field $F$, i.e. $\phi:=h(F,.)$, then a direct calculation shows that the covariant derivative of $F$ can be expressed as follows:
\beqa
h(\tnab_X F, Y)&=&h(X,Y)+\einhalb\Big[d\phi(X,Y) - T^*(F,X,Y)-T^*(F,Y,X)-T^*(X,Y,F)\Big].
\label{nabF}
\eeqa
The homogeneity condition on the torsion implies the following Lemma.

\blem\label{derivdegree}
If $X\in \Gamma^v(T\tem)$ and $Y\in \Gamma^w(T\tem)$, then
$\tnab_XY\in \Gamma^{v+w}(T\tem)$.
\elem
\bprf One derives a Koszul formula for a metric connection with torsion:
\barr{rcrcccl}
2 h(\tnab_X Y, Z)&=&  X\left( h(Y,Z)\right) & +&  Y\left( h(Z,X)\right) & - & Z\left( h(X,Y)\right)\\
&&+\ h\left( [X,Y], Z\right) &-& h\left( [Y,Z], X\right) & + & h\left( [Z,X], Y\right)\\
&&+\ T^*(X,Y,Z)& - & T^*(Y,Z,X) & + & T^*(Z,X,Y).
\earr
By the previous lemma the statement follows.
\eprf

\section{The tractor bundle defined by ambient objects}
In the first part of this section we want to recall the definition of the tractor bundle in terms of an ambient metric due to \cite{cap-gover03}. We do no longer distinguish between $\cal Q$ and the embedded $\iota(\cal Q)$ and between $\vf$ and $\wt{\vf}$.  One considers the restriction on the tangent bundle of an ambient manifold to $\cal Q$, $T\tem |_{\cal Q}$. There is a natural $\rr^+$-action $\vf^*$ on this bundle defined by
\beqa
\vf^*(t, X_q)&:=& t^{-1} \left(d\vf(t,.)\right)_q (X_q).\eeqa
The natural projection $\pi: T\tem |_{\cal Q}\rightarrow \cal Q$ is equivariant w.r.t. this action, i.e.
\[ \begin{CD}
\rr^+ \times T\tem |_{\cal Q} @>{\vf^*}>> T\tem |_{\cal Q} \\
@V{id\times\pi}VV  @VV{\pi}V\\
\rr^+\times {\cal Q} @>{\vf}>> {\cal Q}
\end{CD}\]
commutes. Therefore the quotient 
\be \cal T&:=& \left( T\tem |_{\cal Q}\right)/\rr^+ \ee
is a vector bundle over $M\simeq \cal Q/\rr^+$ of fibre dimension $n+2$. 
In order to induce a metric and a covariant derivative on $\cal T$ by the ambient metric and connection we need the identifications of the next lemma.
\blem There exist the following bijections:\label{identlem}
\bnum
\item $\cal T_{[g]}\simeq \Gamma^{-1}\left( T\tem |_{\rr^+\cdot g}\right)$ for any $[g]\in M=\cal Q/\rr^+$ and $ \rr^+ \cdot g$ the orbit of $g$ under $\vf$;
\item $\Gamma (\cal T)\simeq\Gamma^{-1}\left(T\tem |_{\cal Q}\right)$;
\item For $M=\cal Q/\rr^+$ it is $T_{[q]} M\ \simeq\  \Gamma^0(T\cal Q |_{\rr^+\cdot g}) \text{ mod } \rr\cdot F$;
\item $\Gamma\left(T \left(\cal Q /\rr^+\right)\right)\ \simeq\  \Gamma^0(T\tem |_{\cal Q}) \text{ mod } \Gamma^0(\rr)\cdot F$, where $\Gamma^0(\rr) $ denotes functions which are constant along the orbits of $\rr^+$.
\enum\elem
\bprf
For $g\in \cal Q$ we consider  the equation 
\[ \left(\cal L_F X\right) (\vf_t(g))\ =\ -X(\vf_t(g))\]
along $\rr^+\cdot g=\{\vf_t(g) \mid t\in \rr^+\}$. This
is an ODE which has a unique solution $X(\vf_t(g))=\frac{1}{t}d\vf_t(X_g)$ for each initial value $X_g\in T_g \tem$. Recalling that the equivalence relation in $\cal T$ is the following $X_{\vf_t(g)}=\frac{1}{t} d\vf_t (X_g)$ ensures that there is only one representative lying in $T_g \tem$ and on the other hand that the solutions obtained by starting at different points in $\rr^+\cdot g$ are the same. This proves the first point, with the second immediately following. The third point is obtained in the same manner, and  implies the last one.
\eprf
This gives the following conclusions.
\blem
The ambient metric defines a metric on the bundle $\cal T$.
\elem
\bprf 
For $X\in \cal T_{[g]}$ we denote by $\wt{X}$ the corresponding element in $\Gamma^{-1}\left( T\tem |_{\rr^+\cdot g}\right)$. Lemma \ref{tensordegrees} ensures that $h(X,Y)$ is constant along the orbits of $\rr^+$. Hence, by
\[h^{\cal T}_{[g]}(\wt{X},\wt{Y})\ :=\ h_g(\wt{X}, \wt{Y})\]
a metric on $\cal T$ is defined.
\eprf
\bs
If $d\phi |_{\cal Q}=0$ and $T^*|_{\cal Q}=0$ the ambient connection defines a connection $\nabt$ on  $\cal T$ which parallelises $h^{\cal T}$ and  $(\cal T, h^\cal T, \nabt)$ is a conformal standard tractor bundle for the conformal structure given by $[g]$.
\es
\bprf
By assuming $d\phi |_{\cal Q}=0$ and $T^*|_{\cal Q}=0$ equation (\ref{nabF}) implies that
$\tnab_XF=X$. For $X\in \Gamma^{-1}(T\tem |_{\cal Q}) $ this gives on $\cal Q$ that
\barr{rcccl} 0&=&
\tnab_F X- \tnab _XF-[F,X]
&=&
\tnab_F X
\earr
since $T^* |_{\cal Q}=0$.
Furthermore Lemma \ref{derivdegree} ensures that $\tnab_X$ sends $\Gamma^{-1}(T\tem |_\cal Q)$ to itself if $X\in \Gamma^0(T \tem |_\cal Q )$. Denoting again by a tilde the corresponding identifications of 2. and 4. in Lemma \ref{identlem} we define the connection on $\cal T$ by
\beqa\label{tractcon} \wt{\nabt_X Y}:= \tnab_{\wt{X}} \wt{Y}\eeqa
for $X\in \Gamma(TM)$ and $Y\in \Gamma(\cal T)$. By definition this connection is metric w.r.t. $h^{\cal T}$.

The injection $\cal E[-1]\hookrightarrow \cal T$ is given by sending a function $f$ to $[f\cdot F]\in\cal T$ defining the one-dimensional, light-like subbundle $\cal T^1$. The non-degeneracy condition of Theorem \ref{capgoverthm} follows immediately since $\tnab_X (f F)= X(f) F + f\tnab_XF= X(f) F + f X$ along $\cal Q$.
\eprf

\blem
Let $\tem$ be an ambient manifold with ambient metric $ h$ and ambient connection $\tnab $ such that  $d\phi |_{\cal Q}=0$ and $T^*|_{\cal Q}=0$. Then the  the fundamental vector field $F$ and  the curvature $\wt{\cal R}$ of $\tnab$ satisfy:
\bnum
\item $\wt{\cal R}(X,Y)F=0$ for all $X,Y\in T\cal Q$.
\item If in addition $F\inter T=T(.,.,F)=0$, then 
\bnum
\item $h(\wt{\cal R}(F,Z)X,Y)=h(\wt{\cal R}(X,Y)F,Z)$ along $\cal Q$ for all $X,Y,Z\in T\tem|_{\cal Q}$, and 
\item $
h(\wt{\cal R}(X,Y)F,Z)=0
$
along $\cal Q$ for all $Z\in T\cal Q$, and $X,Y\in T\tem|_{\cal Q}$.
\enum
\enum
\elem
\bprf
A direct calculation based on (\ref{nabF}) gives for the ambient curvature $\wt{\cal R}$ of $F$ that
\beqa
h(\wt{\cal R}(X,Y)F,Z)&=&\left(\Lambda_{(1,2)}\tnab d\phi\right)(X,Y,Z)\nonumber
\\
&&  
+ 2 \left(\Lambda_{(1,2)(3,4)}\tnab T^*\right)(X,Y,F,Z) 
+ \left(\Lambda_{(1,4)}\tnab T^*\right)(X,F,Z,Y)\label{curvF}
\\
&&+\text{ contractions and symmetrisations of $d\phi$ and $T^*$},\nonumber
\eeqa
in which $\Lambda_{(i,j)}$ denotes the skew symmetrisation w.r.t. the $(i,j)$-th component.
Since it was assumed that $d\phi |_{\cal Q} =0$ and $T^*|_{\cal Q}=0$ this implies that 
\beqa\label{curF}
h(\wt{\cal R}(X,Y)F,Z)=0\ \text{  if }\  X,Y \in T\cal Q \text{ and } Z\in T\tem |_\cal Q.\eeqa
This implies that $\wt{\cal R}(X,Y)F=0$ if $ X,Y \in T\cal Q $.

In order to prove (2) first
we notice that for $X,Y,Z\in T\tem|_\cal Q$ the assumption
$F\inter T$ implies
\begin{equation}
\label{dert}(\nabla_XT)(F,Y)\ =\ \tnab_X(T(F,Y))\ =\ 0
\end{equation}
and $T(X,Y)\in F^\bot$ implies
\begin{equation}
\label{dertt}\nabla_XT^*(Y,Z,F)\ =\ X(T^*(Y,Z,F))\ =\ 0
\end{equation}
along $\cal Q$.
Since our connection has torsion, the first Bianchi identity does not hold and we do not have in general that $
h(\wt{\cal R}(X,Y)F,Z)=h(\wt{\cal R}(F,Z)X,Y)$. But in our situation  along $\cal Q$ this equality holds:
for a connection with torsion one gets
\begin{multline}\nonumber
\wt{\cal R}(X,Y)Z +\wt{\cal R}(Y,Z)X +\wt{\cal R}(Z,X)Y
\ = \\
- (\tnab_X T )(Y,Z) - (\tnab_Y T )(Z,X) - (\tnab_Z T )(X,Y)\\ 
+ \underbrace{T(X,T(Y,Z)) + T(Y,T(Z,X))+ T(Z,T(X,Y))}_{= 0\text{ along }\cal Q.}
\end{multline}
This implies along ${\cal Q}$ for $X,Y,Z\in T\tem|_{\cal Q}$ and $F$ the fundamental vector field
\begin{multline}\nonumber
\wt{\cal R}(F,X)Y +\wt{\cal R}(X,Y)F +\wt{\cal R}(Y,F)X
\\ =\ - \tnab_F (T (X,Y)) - \tnab_X (T (Y,F)) - \tnab_Y ( T (F,X))
\ =\ 0
\end{multline}
and
\begin{multline}\nonumber
h(\wt{\cal R}(X,Y)Z +\wt{\cal R}(Y,Z)X +\wt{\cal R}(Z,X)Y ,F)\\
=\ 
-X(T^*(Y,Z,F)) -Y(T^*(Z,X,F))-Z(T^*(X,Y,F))
\
=\
0.\end{multline}
Then, analogously to the torsion free case, one proves that  
 $h(\wt{\cal R}(F,Z)X,Y)=h(\wt{\cal R}(X,Y)F,Z)$ along $\cal Q$ which is (a).

For (b) we use that $dd\phi=0$ implies that 
\be
\left(\Lambda_{(1,2)}\tnab d\phi\right)(X,Y,Z)
&=& \einhalb \left(\tnab_Z d\phi\right)(X,Y)
+\ \begin{minipage}[t]{6cm} terms involving contractions and symmetrisations of $d\phi$ and $T^*$.\end{minipage}
\ee
Hence, we obtain along ${\cal Q}$ for $Z\in T\cal Q$ and $X,Y\in T\tem |_{\cal Q}$ that $ \left(\tnab_Z d\phi\right)(X,Y)=0$. 
By  (\ref{dert}), (\ref{dertt}), and (\ref{curvF}) we get:
\be
h(\wt{\cal R}(X,Y)F,Z)&=&0,
\ee
which is (b) of the second point. 
\eprf

\bs Let $(\cal T, h^\cal T, \nabt)$ defined as above.
Suppose that $d\phi |_{\cal Q}=0$, $T^*|_{\cal Q}=0$, $F\inter T=T^*(.,.,F)=0$, and that
 $\wt{Ric}$ denotes the Ricci curvature of the ambient connection $\tnab$. Then $\nabt$ is the normal standard tractor connection if and only if $\iota^*\wt{Ric}=0$.
%
%If $d\phi |_{\cal Q}=0$, $T^*|_{\cal Q}=0$, then the standard tractor connection is normal. 
%% if and only if 
%If in addition  $F\inter T=T^*(.,.,F)=0$, then
%the Ricci curvature of the ambient connection   $\wt{\Ric}$ satisfies  $\iota^*\wt{\Ric}=0$.
\es

\bprf
First we note that equation (\ref{tractcon}) immediately implies that 
\beqa\label{curvid}
\wt{\cal R^\cal T(X,Y)Z}& =& \wt{\cal R}(\wt{X},\wt{Y})\wt{Z}
\eeqa
where $\wt{\cal R}$ is the ambient curvature and the other tildes denote the corresponding identifications of Lemma \ref{identlem}.
This implies on the one hand that for every $q\in Q$ the fibre $\cal T^1_{[q]}=\rr\cdot [F_q]\subset \cal T_{[q]}$ is annihilated by $\cal R^\cal T(X,Y)$ for $X,Y\in T_{[q]}M$. This is 
one 
condition for the tractor connection to be normal.

%%Now, (1) of Lemma \ref{identlem} identifies $\left(\cal T^1_{[q]}\right)^\bot$ with $\Gamma^{-1}(T\cal Q|_{\rr^+\cdot q})$ (!!!) This together with the third identification of Lemma \ref{identlem} shows that (\ref{curF}) is equivalent to
%%\beqa\label{curtF}
%%h^\cal T(\cal R^\cal T(U,V)F,W)=0\ \text{  if }\  U,V\in T_{[q]}M \text{ and } W\in \cal T_{[q]}^1.\eeqa
%%This relation ensures that $\cal R^\cal T$ leaves $\cal T^1_{[q]}$ invariant for every $[q]\in M$. 

In order to analyse the Ricci condition we fix a basis along $\cal Q$ of the form $(F, S_1, \ldots , S_n, Z) $ with $S_i\in \Gamma(T\cal Q)$ such that $h(S_i,S_j)=\delta_{ij}$ and $Z\in \Gamma (T\tem |_\cal Q)$ such that $h(Z,F)=1$ and $h(Z,S_i)=0$. The Ricci curvature for $X,Y\in T\cal Q$ then is given by
\be
\wt{\Ric}(X,Y)&=&\underbrace{h( \wt{\cal R}(Z,X)Y,F)}
_{=-h( \wt{\cal R}(Z,X)F,Y)=0}+\underbrace{h(\wt{\cal R}(F,X)Y,Z)}_{=h( \wt{\cal R}(Y,Z)F,X)=0}+ \sumi h( \wt{\cal R}(S_i,X)Y,S_i)
\ee
by the previous lemma. Using the identification (\ref{curvid}),  the vanishing of the third summand is equivalent to the vanishing of the Ricci contraction of the tractor curvature on $(\cal T^1)^\bot/\cal T^1\simeq TM[-1]$. This gives the statement.
\eprf

%\section{Examples}

\section{An ambient connection with Cotton--York tensor as torsion}

In this section we want to give an example of an ambient connection which realises the normal conformal Tractor holonomy.

We set $\tem = \mathbb{R} \times M \times \mathbb{R}^+$. Let $s$ and $q$ be the coordinates along $\mathbb{R}$ and $\mathbb{R}^+$, with coordinate vector fields $S$ and $Q$, respectively. Throughout, we assume that we have fixed a metric $g$ in the conformal class and with it the corresponding Levi-Civita connection $\nabla$. The embedding of the bundle $\cal Q$ is given by
\barr{rcrcl}
\iota&:& \cal Q&\rightarrow & \tem\\
&&t^2 g_x&\mapsto&(0,x,t),
\earr
i.e. $\cal Q=\{s=0\}$. This is equivariant w.r.t. the $\rr^+$-action on $\cal Q$ given by $\vf_t(g_x)=t^2g_x$ and on $\tem$ given by 
\be\vf_t (s,x,q)&=&(ts,x,tq).\ee
The fundamental vector field for this action is given by $F(s,x,q)=s S+ q Q$.

In order to define the ambient metric we introduce the following endomorphisms.
Denote by  $\mathsf{P}$ the Schouten tensor of the fixed connection $\nabla$, considered (via $g$) as an endomorphism of $TM$. Using the identification of $T_{(s,x,q)}\tem=\rr S\+ T_x M\+ \rr Q$ we define a bundle automorphism $f$ of $T\tem$ by:
\begin{eqnarray*}
f(S)&=&S\\
f(Q)&=&Q\\
f(X) &=& \left( s \mathsf{P} + qId  \right)^{-1}X, 
%\\l &=& \mathsf{P} \circ f = f \circ \mathsf{P}, \\%h &=& Id, \\%g &=& Id.
\end{eqnarray*}
for $X\in T_xM\subset T_{(s,x,q)}\tem$. 
Of course, $f$ is not defined when $-s/q$ is the inverse of an eigen-value of $\mathsf{P}$. However on $\{0 \} \times M \times \{1 \}$, $-s/q = 0$, so $f$ is defined for  for small $s$.
Notice that $f|_{TM}$ commutes with $\ro$ since $f^{-1}$ does.

For $X\in T_x M$ we can understand $f(X)$ as a lift to $\Gamma(T\tem|_{\rr^+\cdot x})$, which we denote by $\wt{X}$. A direct calculation shows that $\wt{X}$ is homogeneous of degree $-1$, i.e.
$\wt{X}\in \Gamma^{-1}(T\tem|_{\rr^+\cdot x})$, or, if $X\in \Gamma(TM)$, then $\wt{X}\in \Gamma^{-1}(T\tem)$ (c.f. Lemma \ref{identlem}). We use this lift to define a metric:
\be
h(S,Q)&=& 1\\
h(\wt{X}, \wt{Y})&=& g(X,Y)\mbox{, for }X,Y\in TM,\ 
 \mbox{ and $0$ otherwise.}
\ee
This metric is by definition homogeneous of degree $2$, and it also satisfies the second condition posed on an ambient metric: $(\iota^*h)_{q^2 g_x}(X,Y)=h_{(0,x,q)}(X,Y)=q^2g(X,Y)$. Hence, the metric $h$ defines an ambient metric on $\tem$. The form dual to $F$ is given by $\phi=sdq+qds$, hence $d\phi=0$ everywhere.

Using the defined lifts we now want to define an ambient connection on $\tem$ by:
\begin{eqnarray*}\tnab_X \widetilde{Y} &=& \widetilde{\nabla_X Y} -  g(X,Y) S -  \mathsf{P}(X,Y) Q, \\ \tnab_X Q &=& \widetilde{X}, \\ \tnab_X S &=& \widetilde{\mathsf{P}(X)}, \\\tnab_Q X &=& \widetilde{  X},\\\tnab_S X &=&\mathsf{P} (\widetilde{X})\ =\ \wt{P(X)},\end{eqnarray*}and all the other $Q$ and $S$ terms being zero.
\blem \label{nabla:h}
$\tnab h=0$.
\elem
\bprf 
In the following we assume $X,Y,Z\in \Gamma(TM)$ and $U,V\in \{Q,S\}$.
The first thing to notice is that 
$(\tnab_Uh)(V,X)=(\tnab_Xh)(U,V)=0$.
Secondly we get that
\barr{rcccl}
(\tnab_Xh)(Q,\wt{Y})&=&X(h(Q,\wt{Y})-h(\wt{X},\wt{Y})+g(X,Y)h(Q,S)&=&0 \text{ ,  and }\\
(\tnab_Xh)(S,\wt{Y})&=&X(h(S,\wt{Y})-h(\wt{P(X)},\wt{Y})+P(X,Y)h(S,Q)&=&0.
\earr
 Furthermore we have that
\be
(\tnab_Z h)(\wt{X},\wt{Y})
& =&
Z(g(X,Y))-h(\wt{\tnab_Z{X}},\wt{Y}) - h(\wt{X},\wt{\tnab_UY})\\
&&+ g(Z,X)h(S,\wt{Y}) + g(Z,Y)h(S,\wt{X}) + \ro(Z,X)h(Q,\wt{Y}) + \ro(Z,Y)h(Q,\wt{X}) 
\\
&=&(\tnab_Z g)(X,Y)\\
&=&0.
\ee
Finally one notices that the lifts $\wt{X}$ are $Q$ and $S$ parallel:
\begin{eqnarray*}\tnab_Q \widetilde{X} &=& \left( \tnab_Q f \right) (X) + f \circ f (X) \\&=& -Id \circ \left( s \mathsf{P} + qId  \right)^{-2} (X) + \left( s \mathsf{P} + qId  \right)^{-2}(X) \\&=& 0, \\\tnab_S \widetilde{X} &=& \left( \tnab_S f \right) (X) + f \circ P \circ f(X)) \\&=& -\mathsf{P} \circ \left( s \mathsf{P} + qId  \right)^{-2} (X) + \mathsf{P} \circ \left( s \mathsf{P} + qId  \right)^{-2}(X) \\&=& 0.\end{eqnarray*}
This implies that 
$
(\tnab_U h)(\wt{X},\wt{Y})
\ =\
U(g(X,Y))-h(\tnab_U\wt{X},\wt{Y}) - h(\wt{X},\tnab_U\wt{Y})
\ = \ 0.
$\eprf
\blem
The torsion $T$ of $\tnab$ is given by  the Cotton--York tensor of $g$, i.e.
\barr{rclcrcl}
T(X,Y)&=&s\ \wt{CY(X,}Y) &\text{ and }& T^*(X,Y,\wt{Z})&=&sCY(X,Y,Z)\earr
 for $X,Y,Z\in TM$, and zero otherwise, where $CY(X,Y)=(\nabla_X\ro)(Y)-(\nabla_Y\ro)(X)$ is the Cotton--York endomorphism, respectively its dualisation. In particular, $(\cal L_FT^*)=2T^*$,  $T|_{\cal Q}=0$  and $F\inter T=T^*(.,.,F)=0$.  
\elem
\bprf
For $X,Y\in TM$ we calculate:
\be
T(S,Q)&=&0,\\
T(X,Q)&=&0,\\
T(X,S)&=& \widetilde{\mathsf{P}(X)}-\mathsf{P}(\widetilde{X})\ =\ 0\  \mbox{, since $f$ commutes with $\mathsf{P}$},
\ee
which implies $F\inter T=0$.
In order to get the term 
$T(X,Y)$ for $X,Y\in TM$ we need to calculate $\tnab_XY$. It is
\be
\tnab_XY &=& \tnab_X f^{-1}\wt{Y}\\
&=&  \tnab_X \left(s\ro(\wt{Y})+q\wt{Y}\right)\\
&=&  s\tnab_X \wt{\ro(Y)}+q\tnab_X\wt{Y}\\
&=&  s\left(\wt{\nabla_X \ro(Y)}-\ro(X,Y)S-\ro(X,\ro(Y))Q \right)
+q\left( \widetilde{\nabla_X Y} -  g(X,Y) S -  \mathsf{P}(X,Y) Q\right). 
\ee
This implies
\be 
T(X,Y)&=&s\left( \wt{\nabla_X \ro(Y)}-\wt{\nabla_Y \ro(X)}\right) + q \left( \widetilde{\nabla_X Y}- \widetilde{\nabla_X Y}\right) - [X,Y]
%\\
%&=&
%s\left( \wt{\nabla_X \ro(Y)}-\wt{\nabla_Y \ro(X)}\right) + q \wt{[X,Y]}- [X,Y] 
\\
&=&
s\left( \wt{\nabla_X \ro(Y)}-\wt{\nabla_Y \ro(X)}\right) + 
q \wt{[X,Y]}- \wt{f^{-1}([X,Y] )}
\\
&=&
s\left( \wt{\nabla_X \ro(Y)}-\wt{\nabla_Y \ro(X)}\right) 
- \wt{s\ro([X,Y])}
\\
&=&s \ \wt{CY(X,}Y)
\ee
This gives immediately the formula for $T^*$, $T|_{\cal Q}=0$ and $T^*(.,.,F)=0$.
Also homogeneity two is easily verified:
\be
(\cal L_FT^*)(X,Y,\wt{Z})&=& F\big(s\ CY(X,Y,Z)\big)\\
&&
-s\Big( CY(\underbrace{[F,X]}_{=0},Y,\wt{Z})+ CY(X,\underbrace{[F,Y]}_{=0},\wt{Z})+ CY(X,Y,\underbrace{[F,\wt{Z}]}_{=-\wt{Z}})\Big)\\
&=&
2s\ CY(X,Y,\wt{Z}).
\ee
\eprf
A straightforward calculation gives the  curvature of $\tnab$.
\blem
Let $W$ be the Weyl tensor and $CY$ the Cotton--York tensor of the metric $g$. For $X,Y,Z\in TM$ the curvature of $\tnab$ is given by 
\be
\wt{\cal R}(X,Y)\wt{Z}&=&\wt{W(X,Y)Z}-CY(X,Y,Z)\cdot Q,\ee
and all other terms not determined by the symmetries of the curvature terms being zero. In particular, $\wt{\cal R}(U,V)F=0$ along $\cal Q$ for all $U,V\in T\tem|_{\cal Q}$.
\elem
By this we obtain the Ricci-curvature of $\tnab$.
\blem
The Ricci curvature of $\tnab$ vanishes along $\cal Q$.
\elem
\bprf
Let $\{E_i\}_{i=1}^n$ be a basis of $T_x M$ orthonormal w.r.t. $g$. Thus $(S,\wt{E}_1, \ldots , \wt{E}_n, Q)$ is a basis of $T_{(s,x,q)}\tem$ in which the Ricci-curvature in $(s,x,q)$ can be written as
\be
\wt{\Ric}(U,V)&=&h(\wt{\cal R}(S,U)V,Q)+h(\wt{\cal R}(Q,U)V,S)+\sumi h(\wt{\cal R}(\wt{E}_i,U)V,\wt{E}_i).
\ee
This implies $\wt{\Ric}(Q,.)=0$ and $\wt{\Ric}(S,S)=0$. Furthermore we obtain
\be
\wt{\Ric}(X,S)&=&\sumi CY(\wt{E}_i,X,E_i)\ \text{  and }\\
 \wt{\Ric}(X,\wt{Y})&=&\sumi g(W(\wt{E}_i,X)Y,E_i).
 \ee
But along $\cal Q$ we have that $\wt{E}_i=\frac{1}{q}E_i$ and thus
$\wt{\Ric}=0$ along $\cal Q$
\eprf

We can summarize:
\btheo
$(\tem, h, \tnab)$ is an ambient manifold with ambient metric $h$ and ambient connection $\tnab$ with torsion $T$, satisfying $d\phi=0$, $T|_{\cal Q}=0$ and $F\inter T=T^*(.,.,F)=0$. The corresponding  conformal standard Tractor connection is normal.
\etheo

But we can show even more.
\bs The holonomy of $\tem$ is generated by paths in the embedded manifold $\{0 \} \times M \times \{1 \}$, i.e.
$Hol_{(0,x,1)}(T\tem,\tnab)=Hol_{(0,x,1)}(T\tem|_{M},\tnab)$. \es \begin{proof}
Assume now that $\phi: [0,1] \to \{0 \} \times M \times \{1 \}$ is a path with tangent field $Y$. We parallel transport any given vector of $TN$ along $\phi$, getting the equation:
\be0 =&& \tnab_Y (X + aQ + bS) \\=&& \left( \nabla_Y X + a Y + b \mathsf{P}(Y) \right)
 - g(Y,X)S - \mathsf{P}(Y,X)Q + Y(a)Q + Y(b)S .
\ee
for some section $X$ of $TM|_{\phi}$, $a$ and $b$ functions on $M$. The lift of this section into the whole of $\tem$  is $\widetilde{X} + a Q + bS$. To prove that the holonomy transform is the same along any lift of $\phi$, it suffices to prove
\be \tnab_Y \left( \widetilde{X} + a Q + bS \right) &=& 0 \\
\tnab_Q \left( \widetilde{X} + a Q + bS \right) &=& 0 \\
\tnab_S \left( \widetilde{X} + a Q + bS \right) &=& 0.
\ee
However this is the case, since
\be\tnab_Y \left(\widetilde{X} + a Q + bS\right) &=& f\left( \nabla_X Y + aY + b \mathsf{P}(Y) \right) \\&& - g(X,Y) S - \mathsf{P}(X,Y) Q + Y(a)Q + Y(b)S \\&=& 0,\ee
and the other two equalities were proved in Lemma \ref{nabla:h}.
\end{proof}
\begin{theo}The affine holonomy of $\tnab$ is the same as the Tractor holonomy of $\nabt$.\end{theo}
\begin{proof}
We identify $T\tem$ with $\cal T$ -- $\cal T$ being split as in (\ref{T:splitting}) by the choice of $\nabla$ -- via
\begin{eqnarray*}\left( \begin{array}{c} 1 \\ 0 \\ 0 \end{array} \right) \ =\  S,& \left( \begin{array}{c} 0 \\ X \\ 0 \end{array} \right) \ = \ X, &\left( \begin{array}{c} 0 \\ 0 \\ 1 \end{array} \right) \ =\  Q.\end{eqnarray*}Then $\nabt \cong \tnab$ along $M \cong \{ 0 \} \times M \times \{ 1 \}$. Then by the previous proposition, the two connections have same holonomy.\end{proof}\btheo
If there is a metric in the conformal class with vanishing Cotton--York tensor, then $(\tem, h)$ is an ambient metric in the  sense of \cite{cap-gover03}. In this case, the normal conformal tractor holonomy is the holonomy of the Levi-Civita connection of the semi-Riemannian manifold $(\tem, h)$.

If there is a metric in the conformal class which is Einstein,  $(\tem,h)$  gives a Ricci-flat ambient metric in the sense of \cite{fefferman/graham85}. In particular, it gives 
the cone construction in the conformal Einstein case with $S\not=0$ (see \cite{leitner04killing} and \cite{armstrong05}), and the degenerate cone in the Ricci-flat case (see \cite{leistner05}).
\etheo\begin{proof}
If $CY=0$ the connection is torsion free and thus the Levi-Civita connection of $h$.
In the Einstein case, it suffices to restrict attention to the sub-manifold\begin{eqnarray*}\mathbb{R}^+ \left( \{-\mu^{-1} \} \times M \times \{1 \} \right),\end{eqnarray*}when $S + \mu Q$ is the preserved vector, to get the Einstein cone construction. In the Ricci-flat case, $f=\frac{1}{q} Id$, and thus $h=2dsdq+ q^2g$ which is the degenerate cone construction. Note that in the latter case the ambient connection is singularity free.
\end{proof}
\bbem
$\tem$ is \emph{not} conformally invariant:
Notice that is $\mathsf{P} \neq 0$ at $x \in M$, it must have an eigenvalue $q \neq 0$ at $x$. Then the connection is singular along\begin{eqnarray*}\mathbb{R}^+ \left( \{ 1  \} \times M \times \{ -q \} \right).\end{eqnarray*}
This shows that a conformally Ricci-flat metric that is not Ricci-flat itself gives rise to a connection with singularities. These singularities are not removable, as the term
\be
\tnab_Q X &=& \widetilde{  X}
\ee
demonstrates. Thus $\tnab$ is not a conformally invariant construction (even though it has the same holonomy as the conformally invariant Tractor connection).
\ebem

\bbem
In \cite{leistner05} we proved algebraically that the normal Tractor holonomy algebra of a conformal class which contains a metric with zero Cotton--York tensor  is a Berger algebra. The present approach gives a geometric demonstration of this proven fact because
the constructed ambient connection has no torsion if we start the construction with a metric with zero Cotton York tensor. This implies that its curvature satisfies the first Bianchi-identity, and hence, its holonomy algebra is a Berger algebra. \ebem

\bbem
Alternatively, we may define an affine connection on $\wt M$ rather crudely by:\begin{eqnarray*}\tnab_X Y &=& \nabla_X Y -  qg(X,Y) S -  q\mathsf{P}(X,Y) Q, \\\tnab_X Q &=& \frac{1}{q} X, \\\tnab_X S &=& \frac{1}{q} \mathsf{P}(X), \\\tnab_Q X &=& \frac{1}{q} X, \\\tnab_S X &=& 0,\end{eqnarray*}
and all the other $S$ and $Q$ terms being zero. Although this gives no ambient connection with the properties required in the previous sections,
 a direct calculation establishes that $\wt\nabla$ also has same holonomy as the tractor connection. Moreover its torsion only involves terms with $S$ and a tangent vector to $M$  --- unlike the connection of the previous example, which generically has torsion between two vector fields of $M$ --- and has no singularities.
\ebem

\section{Tractor holonomy and the ambient connection}

In this last section, we aim to give some general holonomy related results, which might help in the construction of ambient connections whose holonomy is the same as the Tractor holonomy. The ideal would be to thus generate either a conformally invariant construction, or one with a torsion that allows one to put strong conditions on the holonomy.

As always, let $\tem$ be an ambient manifold with ambient metric $h$ and ambient connection $\tnab$. By $Hol_x(T\tem, \tnab)$ we denote the holonomy group of $\tnab$ at $x\in \tem$. For $U\subset \tem$ a submanifold and $x\in U$ we denote  by $Hol_x(T\tem |_{U},\tem)$ the group generated by loops running in $U$ around $x$.
\blem
Let $\tem$ be an ambient manifold with ambient metric $h$ and ambient connection $\tnab$ such that $d\phi|_{\cal Q}=T|_{\cal Q}=0$. 
%and $F\inter T=T^*(.,.,F)=0$. 
Assume furthermore that
the flow of the $\rr^+$-action is geodesic w.r.t. $\tnab$, i.e. $\vf_t(q)$ is a geodesic for each $q\in \cal Q$.  Suppose we are given a global section $\s\in \Gamma(\cal Q\rightarrow M=\cal Q/\rr^+)$. Then for any $q\in \s(M)$ it holds
\be
Hol_q(T\tem|_{\cal Q},\tnab)& =& Hol_q(T\tem|_{\s(M)},\tnab).
\ee
\elem
\bprf
First we want to collect some consequences of the assumption that the action of $\rr^+$ defines geodesics. 
The tangent vector field to a geodesic $\vf_t(q)$ for $q\in \s(M)\subset\cal Q$, denoted by $\Phi$ is related to the fundamental vector field $F$ as follows
\be
F(\vf_t(q))
&=&
 \frac{d}{ds}\left( \vf_{e^s}\circ \vf_{t}(q)\right)_{|s=0}
\ =\ 
 \frac{d}{ds}\left( \vf_{te^s}(q)\right)_{|s=0}\ =\ 
t\cdot \frac{d}{ds}\left( \vf_{s}(q)\right)_{|s=t}\\
& =& t\cdot \Phi(\vf_t(q))
\ =\ t\cdot \Phi(t)
\ee
Since $\tnab_{\Phi(t)}\Phi(t)=0$ this implies that 
\[
(\tnab_F F)(\vf_t(q))\ =\ \tnab_F (t\Phi)\ =\  F(t)\Phi(\vf_t(q))+t\underbrace{\tnab_F\Phi}_{=0}
\ =\ t\Phi(\vf(t)) \ =\ F(\vf_t(q)) \]
Furthermore, if $X\in \Gamma(T\cal Q|_{\s(M)})$ we denote by $\wt{X}$ the lifted vector field
given by
\[\wt{X}(\vf_t(q))\ =\ (d\vf_t)_q(X(q)).\]
Obviously it is $\wt{F}(\vf(q))=F(\vf_t(q))$. More importantly, the vector fields $\wt{X}$ are homogeneous of degree $0$. Since the torsion and $d\phi$ vanishes on $\cal Q$  (\ref{nabF}) implies that
\[0\ =\ [F,\wt{X}]\ =\  \tnab_F\wt{X}-\tnab_{\wt{X}} F\ =\ \tnab_F\wt{X} - \wt{X}.\]
 
Since $\vf_t(q)$ are geodesics it is
\[\tnab_{\wt{X}}\wt{Y}\ =\ (d\vf_t)(\tnab_XY)\ =\ \wt{\tnab_XY}.\]
Now we take a curve $\wt{\gamma}$ in $\cal Q$. Then $\wt{\gamma}$ is given by $\wt{\gamma}(t)=\vf_{f(t)}(\gamma(t))$ for $f$ a parameter transformation and $ \gamma$ a curve in $\s(M)$. Its tangent vector is
 \be \dot{\wt{\gamma}}(t)&=&\dot{f}(t) \Phi(\wt{\gamma}(t))+d\vf_{f(t)}(\dot{\gamma}(t))\\
 &=&
\frac{ \dot{f}(t)}{f(t)} F(\wt{\gamma}(t))+\wt{\dot{\gamma}(t)}
\\ & =& \frac{ \dot{f}(t)}{f(t)} \wt{F}(\gamma(t))+\wt{\dot{\gamma}(t)}.\ee

For the curve  $\gamma$  in $\s(M)$ now we consider the vector field $Y$ which is assumed to be parallel displaced along $\gamma$.  $Y$ can be written as
\[Y(t)\ =\ a(t)F(\gamma(t))+ X(t)\]
with $a:\rr\rightarrow \rr$ and $X$ a vector field of $T\s(M)$ along $\gamma$. Hence we have
\be
0\ =\ \tnab_{\dot{\gamma}(t)}Y(t)&
=& \dot{a}(t) F(\gamma(t)) +a (t)\tnab_{\dot{\gamma}(t)} F+ \tnab_{\dot{\gamma}(t)}X(t)
\\
&=& \dot{a}(t) F(\gamma(t)) +a(t) \dot{\gamma}(t)+ \tnab_{\dot{\gamma}(t)}X(t)
\ee

Now consider the vector field $U(t)$ along $\wt{\gamma}(t)$ given by
\be U(t) &=& 
a(t) \Phi(\wt{\gamma}(t)) + \frac{1}{f(t)} d\vf_{f(t)}(X(t))\\
&=& \frac{a(t)}{f(t)} F(\wt{\gamma}(t)) + \frac{1}{f(t)} \wt{X(t)}.
\ee
This vector field is parallel along $\wt{\gamma}$, because
 \be
 \tnab_{\dot{\wt{\gamma}}(t)}U(t)&
=&
\left(\frac{\dot{a}(t)}{f(t)}- a(t)\frac{\dot{f}(t)}{f^2(t)} \right) F(\wt{\gamma}(t))- \frac{\dot{f}(t)}{f^2(t)}\wt{X(t)}
\\
&&
+\ \frac{a(t)}{f(t)}
\left( 
\frac{\dot{f}(t)}{f(t)}\underbrace{(\tnab_{F}F)}_{=F}(\wt{\gamma(t)}  
+  
  \underbrace{\tnab_{\wt{\dot{\gamma}}(t)} F}_{={\wt{\dot{\gamma}}(t)}}
   \right)
+\frac{1}{f(t)}\left( \frac{\dot{f}(t)}{f(t)}\underbrace{\tnab_{F} \wt{X(t)}}_{= \wt{X(t)}} + \tnab_{\wt{\dot{\gamma}}(t)}  \wt{X(t)}\right)
\\
&=&
\frac{1}{f(t)}
\left( 
\dot{a}(t) F( \wt{\gamma}(t))
+  a(t)  \wt{\dot{\gamma}(t)} + \tnab_{\wt{\dot{\gamma}}(t)}  \wt{X(t)}
\right)
\\
&=&
\frac{1}{f(t)}
\wt{\left(
\dot{a}(t) F(\gamma(t))
+  a(t)  \dot{\gamma}(t) + \tnab_{\dot{\gamma}(t)} X(t)
\right)}
\\&=&0.
\ee
But for a loop $\wt{\gamma}$ around $q\in \s(M)$, i.e. $\wt{\gamma}(0)=\wt{\gamma}(1)=q$ we have that $f(0)=f(1)=1$ and $\gamma(0)=\gamma(1)=q$ is also a loop. But this implies that $Y(1)=U(1)$ and hence the holonomies are the same.
\eprf
This lemma immediately implies:
\bs
Let $\tem$ be an ambient manifold with ambient metric $h$ and ambient connection $\tnab$ such that $d\phi|_{\cal Q}=T_{\cal Q}=0$, i.e. $(\tem, h, \tnab)$ defines a standard conformal tractor bundle.  If
the flow of the $\rr^+$-action is geodesic w.r.t. $\tnab$, then the conformal tractor holonomy $
Hol_{[q]}(\cal T, \nabt)$ is isomorphic to $Hol_q(T\tem|_{\cal Q},\tnab)$.
\es

Getting these groups to be equal to the full $Hol_q(T\tem,\tnab)$ is more tricky, however, and requires extra considerations beyond the scope of this paper.

%\bibliography{ALGBIB,thomas,GEOBIB,SPINBIB,HOLBIB,CONF}

\end{document}